\documentclass[article,12pt]{amsart}

\usepackage{amsfonts}
\usepackage{amsmath}
\usepackage{amssymb}
\usepackage{mathrsfs}
\usepackage{graphicx}
\usepackage{color}
\usepackage{epsfig}
\usepackage{enumerate}

 \usepackage[foot]{amsaddr}

\numberwithin{equation}{section}
\theoremstyle{plain}
\newtheorem{thm}{Theorem}[section]
\newtheorem{rem}{Remark}[section]

\newcommand{\dE}{\mathbb{E}}
\newcommand{\dR}{\mathbb{R}}

\newcommand{\dP}{\mathbb{P}}
\newcommand{\dC}{\mathbb{C}}

\newcommand{\cB}{\mathcal{B}}

\newcommand{\cS}{\mathcal{S}}
\newcommand{\rI}{\mathrm{I}}

\newcommand{\veps}{\varepsilon}

\newcommand{\ind}{\mbox{1}\kern-.25em \mbox{I}}

\def\build#1_#2^#3{\mathrel{\mathop{\kern 0pt#1}\limits_{#2}^{#3}}}

\def\videbox{\mathbin{\vbox{\hrule\hbox{\vrule height1.4ex \kern.6em\vrule height1.4ex}\hrule}}}
\def\demend{\hfill $\videbox$\\}

\usepackage{xcolor}

\newcommand{\Hm}[1]{\leavevmode{\marginpar{\tiny%
$\hbox to 0mm{\hspace*{-0mm}$\leftarrow$\hss}%
\vcenter{\vrule depth 0.1mm height 0.1mm width \the\marginparwidth}%
\hbox to 0mm{\hss$\rightarrow$\hspace*{-0mm}}$\\\relax\raggedright
#1}}}

\begin{document}

\title[Descents and major index in a random permutation]
{Sharp large deviations for the number of descents and the major index in a random permutation}
\author{Bernard Bercu}
\thanks{The corresponding author is Bernard Bercu, email address: bernard.bercu@math.u-bordeaux.fr}
\address{Universit\'e de Bordeaux, Institut de Math\'ematiques de Bordeaux,
UMR CNRS 5251, 351 Cours de la Lib\'eration, 33405 Talence cedex, France.}
\email{bernard.bercu@math.u-bordeaux.fr}
\author{Michel Bonnefont}
\email{michel.bonnefont@math.u-bordeaux.fr}
\author{Luis Fredes}
\email{luis.fredes@math.u-bordeaux.fr}
\author{Adrien Richou}
\email{adrien.richou@math.u-bordeaux.fr}
\date{\today}

\begin{abstract}
The aim of this paper is to improve the large deviation principle for the number of descents in a random permutation by establishing a sharp large deviation principle of any order.
We shall also prove a sharp large deviation principle of any order for the major index 
in a random permutation.
\end{abstract}


\maketitle

\vspace{1ex}



\maketitle

\vspace{1ex}


\section{Introduction}
\label{S-I}


Let $\cS_n$ be the symmetric group of permutations on the set of integers $\{1,\ldots,n\}$. 
A permutation $\pi_{{ n}} \in \cS_n$ is said to have a descent at position $k \in \{1,\ldots,n-1\}$ if  $\pi_{ n}(k)>\pi_{ n}(k+1)$. 
Denote by $D_n=D_n(\pi_{ n})$ and $M_n=M_n(\pi_{ n})$ the random variables counting the number of descents and the major index of a permutation $\pi_{ n}$ chosen uniformly at random 
over the set $\cS_n$. They are precisely defined, for all $n \geq 2$, by
\begin{equation}
\label{DEFDNMN}
D_{n} = \sum_{k=1}^{n-1} \rI_{\{\pi_{ n}(k)>\pi_{ n}(k+1)\}}
\hspace{1cm}\text{and}\hspace{1cm}
M_{n} = \sum_{k=1}^{n-1} k \rI_{\{\pi_{ n}(k)>\pi_{ n}(k+1)\}}.
\end{equation}
On the one hand, it follows from a beautiful result from Tanny \cite{Tanny1973}
that $D_n$ shares the same distribution as the integer part of the Irwin-Hall distribution.
More precisely, the distribution of $D_n$ is nothing else than the one of the integer part of 
\begin{equation}
\label{DEFSN}
S_n= \sum_{k=1}^n U_k
\end{equation}
where $(U_n)$ is a sequence of independent and identically distributed random variables
with uniform distribution on the interval $[0,1]$. On the other hand, it is well-known 
from an old result due to MacMahon \cite{MacMahon1917} that $M_n$ shares the same distribution as the number of inversions of the permutation $\pi_{ n}$. Consequently, the distribution of
$M_n$ coincides with that of
\begin{equation}
\label{DEFSIGMAN}
\Sigma_n= \sum_{k=1}^n V_k
\end{equation}
where $(V_n)$ is a sequence of independent random variables such that for all $n \geq 1$,
$V_n$ is uniformly distributed on the integers $\{0, \ldots, n -1\}$. It is easy to see from
\eqref{DEFSN} that for all $n \geq 2$,
$$
\dE[D_n]= \frac{n-1}{2} \hspace{1cm}\text{and}\hspace{1cm}
\text{Var}(D_n)=\frac{n+1}{12}.
$$
Moreover, we also have from
\eqref{DEFSIGMAN} that for all $n \geq 2$,
$$
\dE[M_n]= \frac{n(n-1)}{4} \hspace{1cm}\text{and}\hspace{1cm}
\text{Var}(M_n)=\frac{n(n-1)(2n+5)}{72}.
$$
Using a totally different approach involving the number of leaves in a random recursive tree, it 
follows from Bryc, Minda and Sethuraman \cite{Bryc2009} that
the sequence $(D_n/n)$ satisfies a large deviation principle (LDP) with good rate function 
\begin{equation}
\label{DEFID}
I_D(x)=\sup_{ t \in \dR} \bigl\{xt - L_D(t) \bigr\}
\end{equation}
where the cumulant generating function $L_D$ is given by
\begin{equation}\label{DEFLD}
L_D(t)=\log \Bigl(\frac{\exp(t)-1}{t}\Bigr).
\end{equation}
This result was improved in \cite{Bercu2024} by sharp large deviations and concentration
inequalities for the sequence $(D_n/n)$. More recently, M\'eliot and Nikeghbali \cite{Meliot2024}
have shown that the sequence $(M_n/n^2)$ satisfies an LDP with good rate function 
\begin{equation}
\label{DEFIM}
I_M(x)=\sup_{ t \in \dR} \bigl\{xt - L_M(t) \bigr\}
\end{equation}
where, as expected, $L_M$ is the integral of $L_D$,
\begin{equation}
\label{DEFLM}
L_M(t)=\int_{0}^1 L_D(xt) dx.
\end{equation}
The aim of this paper is to go deeper in the analysis of \cite{Bercu2024} and \cite{Meliot2024}
by proving sharp large deviation principle (SLDP) of any order for the sequences $(D_n)$ and 
$(M_n)$ by use of the decompositions \eqref{DEFSN} and \eqref{DEFSIGMAN}. One can observe that a SLDP of order zero for $(M_n)$ was previously established in \cite{Meliot2024}.
\ \\
The paper is organized as follows. Section \ref{S-MR} is devoted our main results
where we establish a SLDP of any order for the sequences $(D_n)$ and 
$(M_n)$. The proof of the SLDP for the sequence $(D_n)$ is given in
Section \ref{S-SLDPDN} while the proof of the SLDP for the sequence $(M_n)$
is carry out in Section \ref{S-SLDPMN}.


\section{Main results}
\label{S-MR}


\noindent
Our first result improves the LDP for the sequence $(D_n/n)$ previously established by Theorem 1.1 in \cite{Bryc2009} as well as Theorem 3.1 in \cite{Bercu2024}. For any positive real number $x$, denote $\{x\}=\lceil x \rceil - x$.
\begin{thm}
\label{T-SLDPDN}
The number of descents $(D_n)$ satisfies a SLDP of any order. More precisely,
for any $x$ in $]1/2,1[$ and for any integer $p\geq 1$, there exists a sequence
$d_{n,1}(x), \ldots,d_{n,p}(x)$ such that for $n$ large enough,
\begin{equation}
\label{SLDPRDN}
\dP \Big(\frac{D_n}{n} \geq x \Big) = \frac{\exp(-nI_D(x)-\{nx\}t_x)}{\sigma_x t_x \sqrt{2\pi n}} 
\left[ 1+\sum_{k=1}^p \frac{d_{n,k}(x)}{n^k} + O\Big(\frac{1}{n^{p+1}}\Big) \right]
\end{equation} 
where the value $t_x$ is the unique solution of $L_D^\prime(t_x)=x$, 
$\sigma_x^2=L_D^{\prime \prime}(t_x)$ and the coefficients $d_{n,1}(x), \ldots,d_{n,p}(x)$ are bounded in $n$ and
may be explicitly calculated as functions of the derivatives of $L_D$ at point $t_x$. 
For example, the first coefficient $d_{n,1}(x)$ is given by
$$
d_{n,1}(x)=\frac{1}{\sigma_x^2}\left(-\frac{1}{t_x^2} -\frac{\{nx\}}{t_x} -\frac{\{nx\}^2}{2}
-\frac{\{nx\}\ell_D(3)}{6 \sigma_x^2}  -\frac{\ell_D(3)}{2  \sigma_x^2 t_x}
+\frac{\ell_D(4)}{8 \sigma_x^2} - \frac{5\ell_D^2(3)}{24 \sigma_x^4}
 \right)
$$
where $\ell_D(3)$ and $\ell_D(4)$ stand for the higher-order derivatives $\ell_D(3)=L_D^{(3)}(t_x)$ 
and $\ell_D(4)=L_D^{(4)}(t_x)$.
\end{thm}


\begin{rem}
The asymptotic expansion \eqref{SLDPRDN} is inspired by the seminal work of Bahadur
and Rao \cite{Bahadur1960} concerning the SLDP for the sample mean, see also
\cite{Bercu2000} for the SLDP for quadratic forms of Gaussian stationary processes. 
\end{rem}

\noindent
Our second result goes deeper than the recent analysis of M\'eliot and Nikeghbali
\cite{Meliot2024} on the large deviation properties of the major index of random standard tableaux. 
It is more complicated than our previous result as it involves the additional prefactor function
\begin{equation}
\label{DEFHMAJ}
H(t)=\frac{1}{2} \big(L_D(t) - t \big).
\end{equation}

\begin{thm}
\label{T-SLDPMAJ}
The major index $(M_n)$ satisfies a SLDP of any order. More precisely,
for any $x$ in $]1/4,1/2[$ and for any integer $p\geq 1$, there exists a sequence
$m_{n,1}(x), \ldots,m_{n,p}(x)$ such that for $n$ large enough,
\begin{equation}
\label{SLDPRMAJ}
\dP \Big(\frac{M_n}{n^2} \geq x \Big) = \frac{\exp(-nI_M(x)+H(t_x))}{\sigma_x t_x \sqrt{2\pi n}} 
\left[ 1+\sum_{k=1}^p \frac{m_{n,k}(x)}{n^k} + O\Big(\frac{1}{n^{p+1}}\Big) \right]
\end{equation} 
where the value $t_x$ is the unique solution of $L_M^\prime(t_x)=x$, 
$\sigma_x^2=L_M^{\prime \prime}(t_x)$ and the coefficients $m_{n,1}(x), \ldots,m_{n,p}(x)$ 
may be explicitly calculated as functions of the derivatives of $L_D$, $L_M$ and $H$ at point $t_x$. 
For example, the first coefficient $m_{n,1}(x)$ is given by
$$
m_{n,1}(x)=\frac{t_x}{2}-\frac{t_x(1+t_x)}{24} - \{n^2 x\} t_x +\frac{t_x\ell_D(1)}{12}  
$$
$$+ \frac{1}{\sigma_x^2} \left( - \frac{1}{t_x^2} +\frac{h(1)}{t_x}
-\frac{h(2)}{2}  -\frac{h^2(1)}{2}  - 
\frac{\ell_M(3)}{2\sigma_x^2t_x} 
+\frac{\ell_M(4)}{8\sigma_x^2} 
+\frac{h(1)\ell_M(3)}{2\sigma_x^2}
- \frac{5\ell_M^2(3)}{24\sigma_x^4}
\right)
$$
where for $1\leq k \leq 4$, $\ell_D(k)=L_D^{(k)}(t_x)$, $\ell_M(k)=L_M^{(k)}(t_x)$ and
$h(k)=H^{(k)}(t_x)$.
\end{thm}

\begin{rem}
We can even be more precise on the behavior of the coefficients $d_{n,k}(x)$ and $m_{n,k}(x)$ in Theorems \ref{T-SLDPDN} and \ref{T-SLDPMAJ}.
More precisely, it is possible to show that $d_{n,k}(x)$ and $m_{n,k}(x)$ are respectively polynomials on the bounded terms $\{nx\}$ and $\{n^2x\}$,
with coefficients only depending on $t_x$ and on the derivatives of $L_D$, $L_M$ and $H$ at point $t_x$.
\end{rem}

\section{Proof of the SLDP for the number of descents}
\label{S-SLDPDN}

We saw in Section \ref{S-I} that the distribution $D_n$ coincides with that of the
integer part of $S_n$ given by \eqref{DEFSN}. It is quite easy to compute the Laplace transform
of $S_n$. As a matter of fact, as $U_1, \ldots, U_n$ are independent and identically distributed 
with uniform distribution on $[0,1]$, we have for all $t \in \dR$,
\begin{equation}
\dE[\exp(tS_n)]= \dE\Big[ \prod_{k=1}^n\exp(tU_k) \Big]= \prod_{k=1}^n \dE[\exp(tU_k)]=\exp\big(nL_D(t)\big)
\label{LAPLACESN}
\end{equation}
where the cumulant generating function $L_D$ is given by 
\eqref{DEFLD}.
Hence, as $I_D(x)= x t_x - L_D(x)$, we have for any $x\in ]1/2, 1[$,
\begin{eqnarray}
    \dP(D_n\geq  nx) &=& \dP( S_n \geq \lceil nx \rceil), \nonumber\\
     &=&  \dE_n \Big[ \exp\Big(\!\!- t_x S_n + nL_D(t_x)\Big ) 
     \ind_{\big\{\frac{S_n}{n} \geq x +\veps_n \big\}} \Big ], \nonumber\\
    &=&\exp(-nI_D(x) ) \dE_n \Big[ \exp\Big( \!\!-nt_x \Big(\frac{S_n}{n} -x \Big)\Big)\ind_{\big\{\frac{S_n}{n} \geq x+\veps_n \big\}} \Big ],
\label{CP}
\end{eqnarray}
where $\dE_n$ stands the expectation under the new probability $\dP_n$ given by
\begin{equation}
\label{NPROBA}
 \frac{d\dP_n}{d\dP}= \exp\big(t_x S_n -n L_D(t_x)\big)
\end{equation}
and $\veps_n= \{nx\}/n$. Denote by $V_n$ the standardized random variables under 
the new probability $\dP_n$,
$$
V_n= \frac{\sqrt{n}}{\sigma_x}  \Big( \frac{S_n}{n}-x \Big).
$$
Moreover, let $f_n$ and $\Phi_n$ be the probability density function and the 
characteristic function of $V_n$ under $\dP_n$, respectively. 
One can easily see that for all $v \in \dR$,
$$f_n (v)=\sqrt{n}\sigma_x \exp \left( \sqrt{n}t_x\sigma_x v +nxt_x-nL_D(t_x)\right) g_n(\sqrt{n}\sigma_x v +nx )$$
where $g_n$ is the probabilty density function of $S_n$ under $\dP$, that is to say the density of the Irwin-Hall distribution which is bounded with compact support. Thus, $f_n$ is clearly square integrable on $\dR$. Accordingly, it follows from Parseval's identity that
\begin{align}
\dE_n \Big[\! \exp\!\Big( \!\!-nt_x \Big(\frac{S_n}{n} -x \Big)\!\Big)\ind_{\big\{\frac{S_n}{n} \geq x+\veps_n \big\}} \Big ]
&= \dE_n \Big[ \exp\big(\!-t_x\sigma_x \sqrt n V_n\big) \ind_{\big\{V_n \geq \frac{\sqrt n \veps_n}{\sigma_x}\big\} } \Big], \nonumber\\
 &= \int_{\dR} \exp\big(\!-\sigma_x t_x \sqrt{n} v\big) \ind_{\big\{ v \geq \frac{\sqrt n \veps_n}{\sigma_x}\big\} }
         f_{n}(v) dv, \nonumber\\
&= \frac{1}{2\pi} \int_{\dR} \frac{\exp\big(\!-(\sigma_x t_x\sqrt n+iv) \frac{\sqrt n \veps_n}{\sigma_x} \big)} {\sigma_x t_x \sqrt n+iv} \Phi_n(v) dv, \nonumber\\
&= \frac{\exp(-t_x\{nx\})}{2\pi\sigma_x t_x\sqrt{n}} \int_{\dR} \Big(\frac{\exp\big(- \frac{ i\{nx\} v}{\sigma_x \sqrt n} \big)} {1+\frac{iv}{\sigma_x t_x \sqrt n}} \Big)\Phi_n(v) dv.
\label{PARSEVAL}
\end{align}
Therefore, we deduce from \eqref{CP} and \eqref{PARSEVAL} that for any $x\in ]1/2, 1[$,
\begin{equation}
\label{MAINDECO}
\dP(D_n\geq  nx)=\frac{\exp(-nI_D(x)-t_x\{nx\})}{2\pi\sigma_x t_x\sqrt{n}} P_n(x)
\end{equation}
where
\begin{equation}
\label{DEFPND}
P_n(x)=\int_{\dR} \Big(\frac{\exp\big(- \frac{ i\{nx\} v}{\sigma_x \sqrt n} \big)} {1+\frac{iv}{\sigma_x t_x \sqrt n}} \Big)\Phi_n(v) dv.
\end{equation}
Hence, we obtain from \eqref{LAPLACESN} and \eqref{NPROBA} that
\begin{eqnarray}
\Phi_n(v)&=& 
\dE_n \Big[ \exp\Big( \frac{iv \sqrt{n}}{\sigma_x}  \Big( \frac{S_n}{n}-x \Big)\Big], \nonumber\\ 
 &=&\dE \Big[ \exp\Big( \frac{iv S_n}{\sigma_x \sqrt n} - \frac{iv\sqrt n x }{\sigma_x} +t_x S_n -nL_D(t_x)\Big) \Big],\nonumber\\
          &=&\exp\Big( n \Big( L_D\big(t_x+ \frac{iv}{\sigma_x \sqrt n}\big) -L_D(t_x)\Big) - \frac{i v \sqrt nx}{ \sigma_x } \Big).
\label{DECOPHIND}
 \end{eqnarray}
Furthermore, the function $L_D$ is infinitely differentiable in a complex neighbourhood of any point $t_x>0$. Then, it follows from Taylor's theorem in the complex plane that for any integer
$p\geq 1$ and for all $v \in \dR$ such that $|v|< \sigma_x \sqrt{n}$,
\begin{equation*}
\log \Phi_n(v)= -\frac{iv \sqrt n x}{\sigma_x}+n\sum_{k=1}^{2p+3}\! \Big(\frac{iv}{\sigma_x \sqrt n}\Big)^k \frac{\ell_D(k)}{k!} + R_{n,p}(v)
\end{equation*}
where for all integer $k\geq 1$, $\ell_D(k)=L_D^{(k)}(t_x)$, and $R_{n,p}(v)$ is a complex remainder term. Since $L_D^\prime(t_x)=x$ and $L_D^{\prime \prime}(t_x)=\sigma^2_x$, the above Taylor expansion reduces to 
\begin{equation}
\label{DECOPHINDS}
\log \Phi_n(v)= -\frac{v^2}{2}+n\sum_{k=3}^{2p+3}\! \Big(\frac{iv}{\sigma_x \sqrt n}\Big)^k \frac{\ell_D(k)}{k!} + 
R_{n,p}(v).
\end{equation} 
Hereafter, one can easily see that for all $t\neq 0$ and for all integer $k \geq 1$,
\begin{equation}
\label{DERIVLD}
L_D^{(k)}(t)=\sum_{\ell=1}^k \frac{a_{k,\ell}}{(e^t-1)^\ell}+ \frac{b_k}{t^k}
\end{equation}
 where the coefficients $a_{k,\ell}$ are the positive integers 
given, for all  $1\leq \ell \leq k$, by
$$
a_{k,\ell}=\frac{1}{\ell}\sum_{i=1}^\ell (-1)^{\ell -i} \binom{\ell}{i} i^k
$$
and $b_k=(-1)^k (k-1)!$. In addition, for all $t>0$ and for all $v \in \dR$,
$|t+iv| \geq t$ and $|\exp(t+iv)-1| \geq \exp(t)-1$, which implies 
from \eqref{DERIVLD} that  for all $k \geq 1$,
\begin{equation}
\label{UNIFUBD}
\left| L_D^{(k)}\Big(t_x +iv\Big) \right|
\leq \sum_{\ell=1}^{k} \frac{a_{k,\ell}}{(e^{t_x}-1)^\ell}+ \frac{|b_{k}|}{t_x^{k}}.
\end{equation}
It follows from \eqref{UNIFUBD} that for all $v \in \dR$ such that $|v|< \sigma_x \sqrt{n}$, 
\begin{equation}
\label{UNIFBREM}
\big| R_{n,p}(v) \big| \leq \frac{M_{x,p} v^{2p+4}}{n^{p+1}}
\end{equation}
for some positive constant $M_{x,p}$ that does not depend on $v$. Denote
$s_n=\varepsilon_n n^{1/6}$ where $\varepsilon_n=n^{-\alpha}$ with $0<\alpha<1/6$ and let
$$
u_n=n\sum_{k=3}^{2p+3}\! \Big(\frac{iv}{\sigma_x \sqrt n}\Big)^k \frac{\ell_D(k)}{k!}  
+ O\left(\frac{v^{2p+4}}{n^{p+1}}\right).
$$
Since $u_n$ is of order $v^3/\sqrt n$, we deduce from the standard Taylor expansion
$$
\exp(u)= 1 + u + \dots + \frac{u^{2p+1}}{(2p+1)!} + O(u^{2p+2})
$$
that for all $v \in \dR$ such that $|v|< s_n$,
\begin{equation}\label{PHIN}
\Phi_n(v)
= \exp\Big(-\frac{v^2}{2}\Big) \left[1+ \sum_{k=1}^{2p+1}
\frac{\psi_{k}(v)}{(\sqrt n)^k} + O\left( \frac{\max({v^{2(p+2)}}, v^{6(p+1)})}{n^{p+1}}\right) \right]
\end{equation}
where the $\psi_{k}$ are polynomials of degree $3k$ and valuation $k+2$, with only  odd powers  for $k$ odd and only  even powers  for $k$ even.
Hence, we obtain from \eqref{PHIN} 
that for all $v \in \dR$ such that $|v|< s_n$,
\begin{equation}
\label{TAYLORPHIND}
\frac{\Phi_n(v)} {1+\frac{iv}{\sigma_x t_x \sqrt n}} 
= \exp\Big(-\frac{v^2}{2}\Big) \left[1+ \sum_{k=1}^{2p+1}
\frac{\varphi_{k}(v)}{(\sqrt n)^k} +  O\left( \frac{\max({v^{2(p+1)}}, v^{6(p+1)})}{n^{p+1}}\right) \right]
\end{equation}
where the $\varphi_{k}$ are also  polynomials of degree $3k$ and valuation $k+2$, with only  odd powers  for $k$ odd and only  even powers  for $k$ even.
For example, 
\begin{align*}
\varphi_{1}(v)&=- \frac{iv}{\sigma_x t_x}- \frac{iv^3\ell_D(3)}{6\sigma_x^3}, \\
\varphi_{2}(v)&= - \frac{v^2}{\sigma_x^2 t_x^2}- \frac{v^4\ell_D(3)}{6\sigma_x^4 t_x}
+\frac{v^4\ell_D(4)}{24\sigma_x^4} - \frac{v^6\ell_D^2(3)}{72\sigma_x^6}.
\end{align*}
Note that in some special regimes of $v$, it may appears that some terms in  the expansions \eqref{PHIN} and \eqref{TAYLORPHIND} may be removed and added in one of the $O$ terms.
Note also that $\varphi_k(v)$ is purely imaginar for $k$ odd and real for $k$ even.
From now on, the integral $P_n(x)$, given by \eqref{DEFPND}, can be
separated into three parts, $P_n(x)=A_n(x)+B_n(x)+C_n(x)$ where
\begin{eqnarray*}
A_n(x)&=&\int_{|v|< s_n} \left(\frac{\exp\big(- \frac{ i\{nx\} v}{\sigma_x \sqrt n} \big)} {1+\frac{iv}{\sigma_x t_x \sqrt n}} \right)\Phi_n(v) dv, \\
B_n(x)&=& \int_{s_n\leq |v|\leq 2\sigma_x \sqrt{n}}  \left(\frac{\exp\big(- \frac{ i\{nx\} v}{\sigma_x \sqrt n} \big)} {1+\frac{iv}{\sigma_x t_x \sqrt n}} \right)\Phi_n(v) dv, \\
C_n(x)&=& \int_{|v|\geq 2\sigma_x \sqrt{n}}  \left(\frac{\exp\big(- \frac{ i\{nx\} v}{\sigma_x \sqrt n} \big)} {1+\frac{iv}{\sigma_x t_x \sqrt n}} \right)\Phi_n(v) dv.
\end{eqnarray*}
We are going to show that the lower part $A_n(x)$ of $P_n(x)$ will complete the asymptotic expansion \eqref{SLDPRDN} via \eqref{MAINDECO}, while the middle part $B_n(x)$ and the upper part $C_n(x)$ of $P_n(x)$ will play a negligible role.
It follows from straightforward calculation previously made in section 4 of \cite{Bercu2024} that for all $t>0$ and for all $v \in \dR$,
\begin{equation}
\label{MAJPHID0}
\Big| \frac{e^{t+iv} -1}{t+iv}\Big|^2=\Big(\frac{e^t -1}{t}\Big)^2
\Big(\frac{t^2}{t^2+v^2}\Big)\Big(1+\frac{2(1-\cos(v)) e^t}{(e^t-1)^2}\Big).
\end{equation}
As $2(1-\cos(v))\leq \min(v^2,4)$, \eqref{MAJPHID0} implies that
\begin{equation}
\label{MAJPHID1}
\Big| \frac{e^{t+iv} -1}{t+iv}\Big|^2
\leq 
\Big(\frac{e^t -1}{t}\Big)^2
\Big(1-\Big(\frac{(e^t-1)^2 - t^2e^t }{(e^t-1)^2} \Big)  
\frac{v^2}{t^2+ v^2}\Big),
\end{equation}
as well as
\begin{equation}
\label{MAJPHID2}
\Big| \frac{e^{t+iv} -1}{t+iv}\Big|^2
\leq
\Big(\frac{e^t -1}{t}\Big)^2
\Big(\frac{t^2}{t^2+v^2}\Big)\Big(1+\frac{4e^t}{(e^t-1)^2}\Big).
\end{equation}
On the one hand, one can observe that 
$$\sigma^2_x=L_D^{\prime \prime}(t_x)=\frac{(e^{t_x}-1)^2 -t_x^2 e^{t_x}}{t_x^2(e^{t_x}-1)^2}.
$$
Hence, we obtain from \eqref{DECOPHIND} and \eqref{MAJPHID1} that, as soon as
$|v|\leq 2\sigma_x \sqrt{n}$,
\begin{equation}
\label{MAJPHID3}
|\Phi_n(v)| \leq \left(1-\Big(\frac{\sigma_x^2t_x^2}{t_x^2+4\sigma_x^2}\Big)\frac{v^2}{n}\right)^{n/2}
\leq \exp\Big(\!-B_x \frac{v^2}{2}\Big)
\end{equation}
where
$$
B_x=\frac{\sigma_x^2t_x^2}{t_x^2+4\sigma_x^2}.
$$
Therefore, we get from \eqref{MAJPHID3} that
\begin{eqnarray}
|B_n(x)| &\leq&
2 \int_{s_{n}}^{2\sigma_x \sqrt{n}}  \exp\Big(\!-B_x \frac{v^2}{2}\Big)dv, \nonumber \\
&\leq & \frac{2}{s_{n}} \int_{s_{n}}^{2\sigma_x \sqrt{n}}  v \exp\Big(\!-B_x \frac{v^2}{2}\Big)dv,
\nonumber \\
&\leq &  \frac{2}{s_{n}B_x}  \exp\Big(\!- B_x\frac{ s_{n}^2}{2}\Big).
\label{MAJBND}
\end{eqnarray}
The upper bound in \eqref{MAJBND} clearly implies that  
$B_n(x)$ goes exponentially fast to $0$.
On the other hand, we find from \eqref{DECOPHIND} and \eqref{MAJPHID2} that
\begin{equation}
\label{MAJPHID4}
|\Phi_n(v)| \leq \left(\frac{t_x^2}{t_x^2+\frac{v^2}{\sigma_x^2 n}}\right)^{n/2}
\left(1+\frac{4e^{t_x}}{(e^{t_x}-1)^2}\right)^{n/2}
\end{equation}
which leads, for all $n>2$, to 
\begin{eqnarray}
|C_n(x)| &\leq& \Big(t_x^2+\frac{4t_x^2 e^{t_x}}{(e^{t_x}-1)^2}\Big)^{n/2}
\int_{|v|\geq 2\sigma_x \sqrt{n}} \Big(\frac{1}{t_x^2+\frac{v^2}{\sigma_x^2 n}}\Big)^{n/2} \!\!dv, \nonumber \\
&\leq & 2 \sigma_x \sqrt{n} \Big(t_x^2+\frac{4t_x^2 e^{t_x}}{(e^{t_x}-1)^2}\Big)^{n/2}
\int_{v\geq 2} \Big(\frac{1}{t_x^2+v^2}\Big)^{n/2} \!\!dv,
\nonumber \\
&\leq &  \sigma_x \sqrt{n} \Big(t_x^2+\frac{4t_x^2 e^{t_x}}{(e^{t_x}-1)^2}\Big)^{n/2}
\int_{v\geq 2} 2v\Big(\frac{1}{t_x^2+v^2}\Big)^{n/2} \!\!dv,
\nonumber \\
&\leq &  \frac{\sigma_x \sqrt{n}}{(n-2)}(t_x^2+4) (C_x)^{n/2} 
\label{MAJCND}
\end{eqnarray}
where
$$
C_x= \frac{t_x^2(e^{t_x}+1)^2}{(t_x^2+4)(e^{t_x}-1)^2}.
$$
One can easily see that for all $t>0$, $t^2e^t<(e^t -1)^2$, which ensures that for all
$x \in ]1/2,1[$, $0<C_x<1$. Consequently, we obtain from \eqref{MAJCND} that 
$C_n(x)$ also goes exponentially fast to $0$. It only remains to evaluate 
the lower part $A_n(x)$ of $P_n(x)$. We deduce from \eqref{DEFPND} and \eqref{TAYLORPHIND}
that
\begin{equation}
\label{TAYLORAND}
A_n(x)=\int_{|v|< s_n} \!\!\!\!\!\!\!\!\exp\Big(- \frac{ i\{nx\} v}{\sigma_x \sqrt n} -\frac{v^2}{2}\Big) \!\left[1+ \sum_{k=1}^{2p+1}
\frac{\varphi_{k}(v)}{(\sqrt n)^k} + O\left( \frac{1+ v^{6(p+1)}}{n^{p+1}}\right)\right] dv.
\end{equation}
For all $a \in \dR$ and for any integer $m\geq 0$, denote
\begin{eqnarray*}
I_m(a)&=& \int_{\dR}  v^{2m} \exp\Big(\!-iav -\frac{v^2}{2}\Big)   dv, \\
J_m(a)&=& \int_{\dR}  v^{2m+1} \exp\Big(\!-iav -\frac{v^2}{2}\Big)  dv.
\end{eqnarray*}
It follows from standard Gaussian calculation that
\begin{eqnarray}
I_m(a)&=& \exp\Big(\!-\frac{a^2}{2}\Big) \int_{\dR}  v^{2m} \exp\Big(-\frac{1}{2}(v+ia)^2\Big)  dv,
\nonumber \\
&=& \exp\Big(\!-\frac{a^2}{2}\Big) \int_{\dR}  (w-ia)^{2m}\exp\Big(-\frac{w^2}{2}\Big)  dw,
\nonumber \\
&=& \exp\Big(\!-\frac{a^2}{2}\Big)  \sum_{k=0}^{2m} \binom{2m}{k} (-ia)^{2n-k} \int_{\dR} w^k\exp\Big(-\frac{w^2}{2}\Big)  dw,
\nonumber \\
&=& \exp\Big(\!-\frac{a^2}{2}\Big)  \sum_{k=0}^{m} \binom{2m}{2k} (-ia)^{2(m-k)} \int_{\dR} w^{2k}\exp\Big(-\frac{w^2}{2}\Big)  dw
\nonumber \\
&=& \sqrt{2 \pi} \exp\Big(\!-\frac{a^2}{2}\Big)  \sum_{k=0}^{m} \binom{2m}{2k} (-1)^{n-k}a^{2(m-k)} \frac{(2k)!}{2^k k!}
\nonumber \\
&=& \frac{\sqrt{2 \pi} (2m)!}{2^m}\exp\Big(\!-\frac{a^2}{2}\Big)  \sum_{k=0}^{m} 
\frac{(-1)^{k}(2a^{2})^k}{(2k)!(m-k)!}.
\label{DEVINA}
\end{eqnarray}
By the same token,
\begin{equation}
J_m(a)=\frac{\sqrt{2 \pi} (2m+1)!}{2^m}(-ia)\exp\Big(\!-\frac{a^2}{2}\Big)  \sum_{k=0}^{m} 
\frac{(-1)^{k}(2a^{2})^k}{(2k+1)!(m-k)!}.
\label{DEVJNA}
\end{equation}
It can be noticed that if we restrict the integration domains of 
$I_m(a)$ and $J_m(a)$ to the interval $[-s_n,s_n]$, and since $s_n$ goes to infinity,  
we shall obtain the same expressions as in
\eqref{DEVINA} and \eqref{DEVJNA} with additional terms that are exponentially small.
Hereafter, we obtain from \eqref{TAYLORAND} together with \eqref{DEVINA} and
\eqref{DEVJNA} with $a=\{nx\}/\sigma_x \sqrt n$ that for any integer $p\geq 1$, there exists a sequence of reals numbers $c_{n,1}(x), \ldots,c_{n,p}(x)$ such that for $n$ large enough,
\begin{equation}
\label{TAYLORFINAND}
A_n(x)=\sqrt{2 \pi} \exp \Big(- \frac{\{nx\}^2}{2 \sigma_x^2 n}\Big)
\left[ 1+\sum_{k=1}^p \frac{c_{n,k}(x)}{n^k} + O\Big(\frac{1}{n^{p+1}}\Big) \right].
\end{equation}
The coefficients $c_{n,1}(x), \ldots,c_{n,p}(x)$ are bounded in $n$ and may be explicitly calculated using the exact
expression of the polynomials $\varphi_{k}$ in \eqref{TAYLORAND} in conjunction with \eqref{DEVINA} and \eqref{DEVJNA}. Finally, it follows from \eqref{MAINDECO} and
\eqref{TAYLORFINAND} that for any $x\in ]1/2, 1[$,
\begin{eqnarray*}
\dP(D_n\geq  nx)&=&\frac{\exp\Big(\!\!-nI_D(x)-t_x\{nx\}- \frac{\{nx\}^2}{2 \sigma_x^2 n}\Big)}{\sqrt{2\pi}\sigma_x t_x\sqrt{n}}
\left[ 1+\sum_{k=1}^p \frac{c_{n,k}(x)}{n^k} + O\Big(\frac{1}{n^{p+1}}\Big) \right], \\
&=& \frac{\exp(-nI_D(x)-t_x\{nx\})}{\sqrt{2\pi}\sigma_x t_x\sqrt{n}}
\left[ 1+\sum_{k=1}^p \frac{d_{n,k}(x)}{n^k} + O\Big(\frac{1}{n^{p+1}}\Big) \right],
\end{eqnarray*}
which completes the proof of Theorem \ref{T-SLDPDN}.
\demend

\section{Proof of the SLDP for the major index}
\label{S-SLDPMN}

We saw in Section \ref{S-I} that the distribution $I_n$ coincides with that $\Sigma_n$ given by \eqref{DEFSIGMAN}. Let $V_1,\ldots,V_n$ be independent random variables such that for all $1 \leq k \leq n$, $V_k$ is uniformly distributed on the integers $\{0, \ldots, k -1\}$.
We have for all $t \in \dR$,
\begin{eqnarray}
\dE[\exp(t\Sigma_n)] &=& \dE\Big[ \prod_{k=1}^n\exp(tV_k) \Big]= \prod_{k=2}^n 
\left( \frac{1}{k} \sum_{\ell=0}^{k-1} e^{t \ell} \right), \nonumber \\
&=& \exp \left( \sum_{k=2}^n \log \left( \frac{1}{k} \sum_{\ell=0}^{k-1} e^{t \ell} \right)\right), 
\nonumber \\
&=& \exp \left( \sum_{k=2}^n \log \left( \frac{1}{k} \Big(\frac{e^{tk} -1}{e^t -1}\Big)\right)\right),
\nonumber \\
&=& \exp \left( \sum_{k=1}^n \left(\log \left( \frac{e^{tk} -1}{tk}\right)-  \log \left( \frac{e^{t} -1}{t}\right)\right)\right)
\nonumber \\
&=& \exp \left( \sum_{k=1}^n L_D(tk)  - nL_D(t) \right).
\label{LAPSIGMAN}
\end{eqnarray}
Consequently, we obtain from \eqref{LAPSIGMAN} that for all $t \in \dR$,
\begin{equation}
\label{LAPLACESIGMAN}
\dE\Big[\exp\Big(\frac{t\Sigma_n}{n}\Big)\Big] = \exp(nL_n(t))
\end{equation}
where
\begin{equation}
\label{DEFLN}
L_n(t)=\frac{1}{n} \sum_{k=1}^n L_D\Big(\frac{tk}{n}\Big)  - L_D\Big(\frac{t}{n}\Big).
\end{equation}
We have for all $t \in \dR$,
\begin{equation}
\label{LDBERN}
L_D(t)= \sum_{k=1}^\infty \frac{B_k t^k}{k k!}=\frac{t}{2}+ \sum_{k=1}^\infty \frac{B_{2k} t^{2k}}{2k (2k)!}.
\end{equation}
The coefficients $B_{k}$ in \eqref{LDBERN} are the well-known Bernoulli numbers. We have taken only even values since the odd Bernoulli numbers are zero except $B_1=1/2$. 
It follows from the Euler-Maclaurin formula that
\begin{eqnarray}
\sum_{k=1}^n L_D\Big(\frac{tk}{n}\Big)  &=& \int_0^n L_D \Big(\frac{tx}{n}\Big) dx+
\frac{1}{2}(L_D(t)-L_D(0)) + \Delta_{n,p}(t) + R_{n,p}(t), \nonumber \\
&=& n \int_0^1 L_D (tx) dx+
\frac{1}{2}L_D(t)+ \Delta_{n,p}(t) + R_{n,p}(t)
\label{MACLAURIN}
\end{eqnarray}
where
\begin{eqnarray}
\label{DEFDELTANP}
\Delta_{n,p}(t)&=&\sum_{k=1}^p \frac{B_{2k}}{(2k)!} \Big(\frac{t}{n}\Big)^{2k-1} \Big(L_D^{(2k-1)}(t)
-L_D^{(2k-1)}(0)\Big), \nonumber \\
&=& \frac{t}{12n} \Big(L_D^{\prime}(t) - \frac{1}{2} \Big) + \sum_{k=2}^p \frac{B_{2k}}{(2k)!} \Big(\frac{t}{n}\Big)^{2k-1} L_D^{(2k-1)}(t)
\end{eqnarray}
and the remainder term
\begin{equation}
\label{DEFRNPBERN}
R_{n,p}(t)= \frac{1}{(2p+1)!}\Big(\frac{t}{n}\Big)^{2p+1} \int_0^n L_D^{(2p+1)} \Big(\frac{tx}{n}\Big)\cB_{2p+1}(x-\lfloor x \rfloor)  dx
\end{equation}
where $\cB_{p}$ stands for the Bernoulli polynomial of order $p$. We have for any $k \geq 0$,
\begin{equation}
\label{DEFRNP}
|R^{(k)}_{n,p}(t)|= O\Big(\frac{1}{n^{2p}}\Big).
\end{equation}
Hence, we deduce from \eqref{DEFLN}
and \eqref{MACLAURIN} that for all $t \in \dR$, 
\begin{equation}
\label{MAINDECOM}
L_n(t)=L_M(t) +\frac{1}{n} H(t) +\frac{t}{2n}-L_D\Big(\frac{t}{n}\Big)+\frac{1}{n}\Delta_{n,p}(t)
+\frac{1}{n}R_{n,p}(t)
\end{equation}
where the function $L_M$ is given in \eqref{DEFLM} and
the function $H$ is defined in \eqref{DEFHMAJ}.
Hereafter, the stratregy is quite different and more difficult to handle than the one for the descents.
First of all, we have for any $x\in ]1/4, 1/2[$,
\begin{equation}
\label{DECOMAJPR}
\dP(M_n \geq n^2 x) = \dP(\Sigma_n \geq n^2 x) =\sum_{k=\lceil n^2x \rceil }^{s_n} \dP\left(\Sigma_n =k\right)
\end{equation}
where $s_n=n(n-1)/2$. Moreover, we have from \eqref{LAPLACESIGMAN} that for all $t,v\in \dR$,
$$
\exp\big( n L_n(n(t+iv)) \big)=\dE\Big[\exp\big((t+iv) \Sigma_n\big) \Big]= \sum_{k=0}^{s_n} 
\exp\big((t+iv)k\big) \dP(\Sigma_n =k),
$$
which implies that for all $t,v\in \dR$ and for all $k \geq 0$, 
\begin{equation}\label{COEFF-FOURIER}
    \dP\left(\Sigma_n =k\right) =\exp(-tk)  \frac{1}{2\pi} \int_{-\pi}^{\pi} \exp\big( n L_n(n(t+iv)) \big) \exp(-ik v)dv.
\end{equation}
One can observe that the identity \eqref{COEFF-FOURIER} is also true for all $k \geq s_n$.
Hence, we deduce from Fubini's theorem together with \eqref{DECOMAJPR} and
\eqref{COEFF-FOURIER} that for all $t>0$,
\begin{equation}\label{PROBAGEQX}
   \dP(M_n \geq n^2 x)  =    \frac{1}{2\pi} \int_{-\pi}^{\pi} \exp\big( n L_n(n(t+iv)) \big)   \sum_{k = \lceil n^2x \rceil }^{+\infty}   \exp\big(-k(t+i v)\big)dv.
\end{equation}
Accordingly, if we choose $t=t_x/n$, which is positive since $x>1/4$, we find from \eqref{PROBAGEQX} that the probability $\dP(M_n \geq n^2 x)$ can be separeted into two terms,
\begin{equation}
    \dP(M_n\geq  n^2x) = \frac{1}{2\pi} \int_{-\pi}^{\pi} \exp\big( n L_n(t_x+inv) \big)   \Lambda_n(t_x+inv)dv= A_n(x) B_n(x),
\label{CPM}
\end{equation}
where
\begin{equation}
\label{DEFANINI}
A_n(x)=\exp\big( n\big(L_n(t_x) -xt_x\big)\big)
\end{equation}
and
\begin{equation}
\label{DEFBNINI}
B_n(x)= \frac{1}{2\pi} \int_{-\pi}^{\pi} \exp\big( n \big(L_n(t_x+inv) -L_n(t_x)+xt_x\big)\big)   \Lambda_n(t_x+inv)dv
\end{equation}
with
\begin{equation}
\label{DEFLAMBDANPR}
\Lambda_n(t_x+inv)
= \sum_{k = \lceil n^2x \rceil }^{+\infty}   \exp\big(-k(t_x+i nv)/n\big)
=\frac{\exp\big( -\frac{\lceil n^2x \rceil}{n}  (t_x+i nv) \big)}{1-\exp\big(-\frac{1}{n}(t_x+i nv)\big)}.
\end{equation}
Concerning the first term $A_n(x)$, the point $t_x$ is chosen
such that $L_M^{\prime}(t_x)=x$ which ensures that
$I_M(x)=xt_x-L_M(t_x)$, leading via \eqref{DEFANINI} to
\begin{equation}
\label{DEVAN1}
A_n(x)=\exp(-nI_M(x)) \exp\Big( n\big(L_n(t_x) -L_M(t_x)\big)\Big).
\end{equation}
Moreover, we have from \eqref{MAINDECOM} that
\begin{equation}
\label{DEVAN2}
n\big(L_n(t_x) -L_M(t_x)\big)=H(t_x)+\frac{t_x}{2}-nL_D\Big(\frac{t_x}{n}\Big)+\Delta_{n,p}(t_x)
+R_{n,p}(t_x).
\end{equation}
We already saw from \eqref{LDBERN} that
\begin{equation}
\label{DEVAN3}
nL_D\Big(\frac{t_x}{n}\Big)=\frac{t_x}{2}+ n\sum_{k=1}^p \frac{B_{2k}}{2k(2k)!} 
\Big( \frac{t_x}{n}\Big)^{2k}+O\Big(\frac{1}{n^{2p}}\Big). 
\end{equation}
As in Section \ref{S-SLDPDN}, for all integer $k\geq 1$, denote 
$\ell_D(k)=L_D^{(k)}(t_x)$ and $\ell_M(k)=L_M^{(k)}(t_x)$. One needs to be very careful here. We already saw that $t_x$ is chosen such that $L_M^{\prime}(t_x)=x$ and we set $\sigma^2_x=L_M^{\prime \prime}(t_x)$.
This will not provide tractable information on $\ell_D(k)$ as we recall from \eqref{DEFLM} that 
$L_M$ is the integral of $L_D$,
\begin{equation}
L_M(t)=\int_{0}^1 L_D(xt) dx=\frac{t}{4}+ \sum_{k=1}^\infty \frac{B_{2k}t^{2k}}{2k(2k+1)!}.
\end{equation}
However, \eqref{DEFDELTANP} reduces to
\begin{equation}
\label{DEVAN4}
\Delta_{n,p}(t_x)=\sum_{k=1}^p \frac{B_{2k} \ell_D(2k-1)}{(2k)!} \Big(\frac{t_x}{n}\Big)^{2k-1} 
-\frac{t_x}{24 n}.
\end{equation}
Hence, it follows from \eqref{DEVAN1} together with the three contributions \eqref{DEVAN2},  \eqref{DEVAN3} and \eqref{DEVAN4} that
\begin{equation}
\label{DEVANFIN}
A_n(x)=\exp\big(\!-nI_M(x)+H(t_x) \big)
\left[ 1+\sum_{k=1}^p \frac{a_{k}(x)}{n^k} + O\Big(\frac{1}{n^{p+1}}\Big) \right]
\end{equation}
where the prefactor was previously defined in \eqref{DEFHMAJ}
and $a_{1}(x), \ldots,a_{p}(x)$
can be explicitly calculated as functions of the derivatives of $L_D$ at point $t_x$. 
For example, $a_{1}(x)$ and $a_{2}(x)$ are given by
\begin{align*}
   a_{1}(x) &= - \frac{t_x(1+t_x)}{24} + \frac{t_x \ell_D(1)}{12}, \\
   a_{2}(x) &= \frac{a_{1}(x)^2}{2}.
\end{align*}
We shall now focus our attention on the second term $B_n(x)$. By the change of variable $w=n \sqrt{n} v$, we have from \eqref{DEFBNINI} that
$$
B_n(x) = \frac{1}{2\pi n \sqrt{n}} \int_{-\pi n \sqrt{n}}^{\pi n \sqrt{n}} 
\!\!\exp\left(\! n \Big(L_n\Big(t_x+\frac{iw}{\sqrt{n}}\Big) -L_n(t_x)+xt_x\Big)\!\right)\!
\Lambda_n\Big(t_x+\frac{iw}{\sqrt{n}}\Big) dw.
$$
Denote
$s_n=\varepsilon_n n^{1/6}$ where $\varepsilon_n=n^{-\alpha}$ with $0<\alpha<1/6$.
As previously done, $B_n(x)$ can be separated into three terms,
$$
B_n(x)=\frac{1}{ \sqrt{2 \pi n}} \big(C_n(x)+D_n(x)+E_n(x))
$$
where
\begin{align*}
C_n(x) &= \frac{1}{n\sqrt{2 \pi}} \int_{|w| < s_n} \!\!\!\!\exp\left(\! n \Big(L_n\Big(t_x+\frac{iw}{\sqrt{n}}\Big) -L_n(t_x)+xt_x\Big)\!\right)\!
\Lambda_n\Big(t_x+\frac{iw}{\sqrt{n}}\Big) dw, \\
D_n(x) &= \frac{1}{n\sqrt{2 \pi}} \int_{s_n \leq |w| \leq a \sqrt{n}} \!\!\!\!\!\!\!\!\exp\left(\! n \Big(L_n\Big(t_x+\frac{iw}{\sqrt{n}}\Big) -L_n(t_x)+xt_x\Big)\!\right)\!
\Lambda_n\Big(t_x+\frac{iw}{\sqrt{n}}\Big) dw, \\
E_n(x) &= \frac{1}{n\sqrt{2 \pi}} \int_{a \sqrt{n} \leq |w| \leq \pi n \sqrt{n}} \!\!\!\!\!\!\!\!\exp\left(\! n \Big(L_n\Big(t_x+\frac{iw}{\sqrt{n}}\Big) -L_n(t_x)+xt_x\Big)\!\right)\!
\Lambda_n\Big(t_x+\frac{iw}{\sqrt{n}}\Big) dw
\end{align*}
and $a$ is a positive parameter that will be chosen later. 
It is well-known that for all $z \in \dC$ with $z \neq 0$,
$$
\frac{1}{1- \exp(-z)}= \frac{1}{z}\sum_{k=0}^\infty \frac{B_kz^k}{k!}.
$$
Consequently, we obtain from \eqref{DEFLAMBDANPR} that
\begin{align}
C_n(x)
&= \frac{\exp ( -\tau_n(x))}{n\sqrt{2 \pi}} \int_{|w|<s_n} 
\!\!\exp\left(\! n \Big(L_n\Big(t_x+\frac{iw}{\sqrt{n}}\Big) -L_n(t_x)-\frac{iwx}{\sqrt{n}}\Big)\!\right)\nonumber \\
& 
\hspace{0.5cm} 
 \exp \Big(\!- \frac{{i\{n^2 x \}w}}{n \sqrt{n}} \Big)
\left( 1-\exp\Big(\!-\frac{t_x}{n}-\frac{iw}{n\sqrt{n}}\Big)\right) ^{-1} \!\!dw,
\nonumber \\
&= \frac{\exp ( -\tau_n(x))}{ t_x\sqrt{2 \pi}} \int_{|w|<s_n} 
\!\!\exp\left(\! n \Big(L_n\Big(t_x+\frac{iw}{\sqrt{n}}\Big) -L_n(t_x)-\frac{iwx}{\sqrt{n}}\Big)\!\right)\nonumber \\
& 
\hspace{0.5cm}  
\exp \Big(\!- \frac{{i\{n^2 x \}w}}{n \sqrt{n}} \Big)
\left( 1+\frac{iw}{t_x\sqrt{n}}\right) ^{-1} \sum_{k=0}^\infty \frac{B_k t_x^k}{k!n^k} \left(  1+\frac{iw}{t_x\sqrt{n}} \right)^k\!\!dw,
\label{DEVCN1}
\end{align}
where
$$
\tau_n(x)= \frac{1}{n} \{n^2 x \} t_x.
$$
It follows Taylor's theorem that for all $|w| <s_n$,
\begin{equation}
\label{DEVBN1}
L_n\Big(t_x+\frac{iw}{\sqrt{n}}\Big)-L_n(t_x) -\frac{iw}{\sqrt{n}} L_n^\prime (t_x)= \sum_{j=2}^{2p+3}\frac{1}{j!} \Big(\frac{iw}{\sqrt{n}}\Big)^j L_n^{(j)}(t_x) + \zeta_{n,p}(w)
\end{equation}
where the remainder $\zeta_{n,p}(w)$ satisfies for all $w \in \dR$ such that $|w| <s_n$,
\begin{equation}
\label{NEWREM}
| \zeta_{n,p}(w) | \leq \frac{M_{x,p} w^{2p+4}}{n^{p+2}}
\end{equation}
for some constant positive $M_{x,p}$ that does not depend on $w$.
We obtain from \eqref{LDBERN}, \eqref{DEFDELTANP}, \eqref{DEFRNP} and \eqref{MAINDECOM} that
\begin{align*}
L_n^\prime (t_x)&=L^\prime_M(t_x) +\frac{H^\prime(t_x) }{n} +\frac{1}{2n}-\frac{1}{n}L_D^\prime\Big(\frac{t_x}{n}\Big)+\frac{1}{n}\Delta^\prime_{n,p}(t_x)
+\frac{1}{n}R^\prime_{n,p}(t_x), \\
&=x +\frac{H^\prime(t_x) }{n}   -\frac{1}{24n^2} +\frac{1}{n} \sum_{k=1}^p \frac{B_{2k}}{(2k)!} \Big(\frac{t_x}{n}\Big)^{2k-1} \big(\ell_D(2k)-1\big) \\
&\hspace{1cm}+\frac{1}{n^2}\sum_{k=1}^p \frac{(2k-1) B_{2k}}{(2k)!}  \ell_D(2k-1) \Big(\frac{t_x}{n}\Big)^{2k-2} 
\!\!\!\!+O\Big(\frac{1}{n^{2p+1}}\Big),
\end{align*}
which reduces to
\begin{equation}
\label{DEVBNL1}
L_n^\prime (t_x) =x +\frac{H^\prime(t_x) }{n}   -\frac{1}{24n^2} +\sum_{k=1}^p \frac{\psi_{k,1}(x)}{n^{2k}} + O\Big(\frac{1}{n^{2p+1}}\Big),
\end{equation}
where
$$
\psi_{k,1}(x)= \frac{B_{2k} (t_x)^{2k-2}}{2k (2k-2)!} 
\left( \frac{t_x(\ell_D(2k) -1)}{2k-1} +  \ell_D(2k-1) \right).
$$
By the same token,
for any $ j \geq 2$,
\begin{align}
L_n^{(j)} (t_x)&=L^{(j)}_M(t_x) +\frac{H^{(j)}(t_x) }{n} -\frac{1}{n^j}L_D^{(j)}\Big(\frac{t_x}{n}\Big)+\frac{1}{n}\Delta^{(j)}_{n,p}(t_x)
+\frac{1}{n}R^{(j)}_{n,p}(t_x),  \nonumber \\
&=\ell_M(j)+\frac{H^{(j)}(t_x) }{n} +\sum_{k=1}^p \frac{\psi_{k,j}(x)}{n^{2k}} + O\Big(\frac{1}{n^{2p+1}}\Big),
\label{DEVBNLj}
\end{align}
where the coefficients $\psi_{k,j}(x)$ may be explicitly calculated as functions of the derivatives of $L_D$ at point $t_x$. For example, in the special case $j=2$, as $\ell_M(2)=\sigma^2_x$,
\begin{equation}
\label{DEVBNL2}
L_n^{\prime \prime} (t_x) = \sigma^2_x+\frac{H^{\prime \prime}(t_x) }{n} +\sum_{k=1}^p \frac{\psi_{k,2}(x)}{n^{2k}} + O\Big(\frac{1}{n^{2p+1}}\Big),
\end{equation}
where $\psi_{1,2}(x)=\big(t_x \ell_D(3) + 2 \ell_D(2) -1\big)/12$ and for $k \geq 2$,
$$
\psi_{k,2}(x)= \frac{B_{2k} (t_x)^{2k-3} }{2k (2k-3)!} 
\left(\frac{t_x^2\ell_D(2k +1)}{(2k-1)(2k-2)}  + \frac{t_x (2\ell_D(2k)-1)}{(2k-2)} + \ell_D(2k-1)  \right).
$$
Then, we deduce from \eqref{DEVBNL1}, \eqref{DEVBNLj} and \eqref{DEVBNL2} that
\begin{align*}
&n \Big(L_n\Big(t_x+\frac{iw}{\sqrt{n}}\Big)-L_n(t_x) -\frac{iw x}{\sqrt{n}}\Big)= \frac{iw H^\prime(t_x)}{\sqrt{n}} -\frac{iw}{24n\sqrt{n}} -\frac{w^2 \sigma^2_x}{2} 
-\frac{w^2  H^{\prime \prime}(t_x)}{2n} \\
&+\sum_{j=3}^{2p+3} \frac{1}{j!} \Big(\frac{iw}{\sqrt{n}}\Big)^j \left( n\ell_M(j)+H^{(j)}(t_x) \right) 
+\sum_{j=1}^{2p+3} \frac{n}{j!} \Big(\frac{iw}{\sqrt{n}}\Big)^j \left( \sum_{k=1}^p \frac{\psi_{k,j}(x)}{n^{2k}} + O\Big(\frac{1}{n^{2p+1}}\Big)\right) 
\nonumber \\& 
+ n\zeta_{n,p}(w). \nonumber
\end{align*}
Consequently, similarly to \eqref{TAYLORAND},  we obtain from \eqref{DEVCN1} that 
\begin{align}
\label{DEVCN2}
C_n(x) &=  \frac{\exp ( -\tau_n(x))}{t_x\sqrt{2 \pi}} \!\!\int_{|w|<s_n} 
\exp\Big(\!\! -\frac{w^2 \sigma^2_x}{2}  - \frac{{i\{n^2 x \}w}}{n \sqrt{n}} \Big) \notag\\
& \left[1+ \sum_{k=1}^{2p+1}
\frac{\varphi_{k}(w)}{(\sqrt n)^k} + O\left(\frac{1+ w^{6(p+1)}}{n^{p+1}}\right)\right]\! dw
\end{align}
where the $\varphi_{k}(w)$ are polynomials in odd powers of $w$ for $k$ odd, and polynomials 
in even powers of $w$ for $k$ even. For example,
\begin{align*}
\varphi_{1}(w)&= iw H^\prime(t_x)- \frac{iw}{t_x} - \frac{iw^3\ell_M(3)}{6}, \\
\varphi_{2}(w)&= \frac{t_x}{2}  - \frac{w^2}{ t_x^2} -\frac{w^2H^{\prime \prime}(t_x) }{2}  
-\frac{w^2(H^{\prime}(t_x))^2 }{2}  +\frac{w^2 H^\prime(t_x)}{t_x}- \frac{w^4\ell_M(3)}{6 t_x}\\
& +\frac{w^4\ell_M(4)}{24} +\frac{w^4 H^\prime(t_x)\ell_M(3)}{6}
- \frac{w^6\ell_M^2(3)}{72}.
\end{align*}
Note also that $\varphi_k(v)$ is purely imaginar for $k$ odd and real for $k$ even.
Therefore, we deduce from \eqref{DEVCN2}, \eqref{DEVINA} and \eqref{DEVJNA} that
\begin{align*}
C_n(x) &=  \frac{\exp \Big(\!\! -\tau_n(x)- \frac{\{nx\}^2}{2 \sigma_x^2 n^3} \Big)}{\sigma_x t_x} 
\left[1+ \sum_{k=1}^{p} \frac{c_{n,k}(x)}{n^k} + O\Big(\frac{1}{n^{p+1}}\Big)
\right],\\
&= 
\frac{1}{\sigma_x t_x} 
\left[1+ \sum_{k=1}^{p} \frac{d_{n,k}(x)}{n^k} + O\Big(\frac{1}{n^{p+1}}\Big)
\right],
\end{align*}
where the coefficients $c_{n,1}(x), \ldots,c_{n,p}(x)$
may be explicitly calculated using the exact expression of the polynomials $\varphi_{k}$
in \eqref{DEVCN2}. For example, $c_{n,1}(x)$ is given by  
\begin{align*}
c_{n,1}(x)&=\frac{t_x}{2}  - \frac{1}{\sigma_x^2 t_x^2} -\frac{H^{\prime \prime}(t_x) }{2\sigma_x^2}  
-\frac{(H^{\prime}(t_x))^2 }{2\sigma_x^2}  +\frac{H^\prime(t_x)}{\sigma_x^2t_x}- 
\frac{\ell_M(3)}{2\sigma_x^4t_x} \\
&+\frac{\ell_M(4)}{8\sigma_x^4} 
+\frac{H^\prime(t_x)\ell_M(3)}{2\sigma_x^4}
- \frac{5\ell_M^2(3)}{24\sigma_x^6}.
\end{align*}
In order to complete the proof of Theorem \ref{T-SLDPMAJ}, it only remains to show that the two last terms $D_n(x)$ and $E_n(x)$ will play a negligible role.
We start by noticing that
\begin{align}\label{BDn}
    |D_n(x)| &\leq 
    \frac{\chi(x,n)}{n\sqrt{2 \pi}} \int_{s_n \leq |w| \leq a \sqrt{n}} \!\left \vert\exp\left(\! n \Big(L_n\Big(t_x+\frac{iw}{\sqrt{n}}\Big) -L_n(t_x)\Big)\!\right)\right \vert dw,
\end{align}
where we use the inequality $\exp(nxt_x)|\Lambda_n(t_x+iw/\sqrt{n})|\leq \chi(x,n)$ with
\begin{equation}
\label{DEFCHI}
\chi(x,n) = \exp(-\tau_n(x))\left\vert 1-\exp\left(\frac{-t_x}{n}\right)\right\vert^{-1} \leq \frac{n e^{t_x}}{t_x}.
\end{equation}
\\
From now on, we focus on the integral appearing in \eqref{BDn}. We have from \eqref{DEFLN} that
\begin{align*}
L_n\Big(t_x+\frac{iw}{\sqrt{n}}\Big)  -L_n(t_x) =& \frac{1}{n}\sum_{k=1}^n 
\left(L_D\Big(\frac{k}{n}\Big(t_x+\frac{iw}{\sqrt{n}}\Big)\Big)  - L_D\Big(\frac{kt_x}{n}\Big)\right) \\
&-
\left( L_D\Big(\frac{1}{n}\Big(t_x+\frac{iw}{\sqrt{n}}\Big)\Big)-L_D\Big(\frac{t_x}{n}\Big) \right),
\end{align*}
which leads to
\begin{align}
\label{MAJDNXNUM0}
\exp\left(\! n \Big(L_n\Big(t_x+\frac{iw}{\sqrt{n}}\Big) -L_n(t_x)\Big)\!\right)\!
&=\frac{{\displaystyle\prod_{k=1}^n}\exp\left(L_D\Big(\frac{k}{n}\Big(t_x+\frac{iw}{\sqrt{n}}\Big)\Big)  - L_D\Big(\frac{kt_x}{n}\Big)\right) }{\exp\left( n \Big(L_D\Big(\frac{1}{n}\Big(t_x+\frac{iw}{\sqrt{n}}\Big)\Big)-L_D\Big(\frac{t_x}{n}\Big)\Big) \right)}.
\end{align}
Hereafter, our goal is to find a suitable upper bound for the numerator in equation \eqref{MAJDNXNUM0}. It follows from inequality \eqref{MAJPHID1} that as soon as
$|w|\leq a \sqrt{n}$,
\begin{equation}
\label{MAJDNXNUM1}
\left\vert \prod_{k=1}^n  \exp\left (L_D\Big(\frac{k}{n}\Big(t_x+\frac{iw}{\sqrt{n}}\Big)\Big)  - L_D\Big(\frac{kt_x}{n}\Big)\right) \right\vert
 \leq \exp\Big( -\sum_{k=1}^n C\Big(\frac{kt_x}{n}\Big) \frac{w^2}{2n} \Big)
\end{equation}
where for all $t>0$,
$$
C(t)=\frac{t^2}{t^2 + a^2}L_D^{\prime \prime}(t).
$$
One can observe that $L_D^{\prime \prime}(t)$ is a strictly decreasing function. Consequently, we deduce from 
\eqref{MAJDNXNUM1} together with a simple  Maclaurin-Cauchy test that 
\begin{equation}
\label{MAJDNXNUM2}
\left\vert \prod_{k=1}^n  \exp\left (L_D\Big(\frac{k}{n}\Big(t_x+\frac{iw}{\sqrt{n}}\Big)\Big)  - L_D\Big(\frac{kt_x}{n}\Big)\right) \right\vert
 \leq \exp\Big( - D_{a,x} L_D^{\prime \prime}(t_x)\frac{w^2}{2} \Big)
\end{equation}
where $D_{a,x}$ is the positive constant given by
$$
D_{a,x}=\int_0^1 \frac{(t_x y)^2}{(t_x y)^2 + a^2} dy.
$$
Hence, it will be enough to show that the denominator of equation \eqref{MAJDNXNUM0} is bounded away from zero to deduce that $D_n(x)$ is exponentially negligible. We have from Taylor's theorem in the complex plane that
$$
nL_D\left(\frac{z}{n}\right) = \frac{z}{2} + O\left(\frac{z^2}{n}\right).
$$
Therefore, for all $w \in \dR$ such that $|w| \leq a \sqrt{n}$,
$$
n \Big(L_D\Big(\frac{1}{n}\Big(t_x+\frac{iw}{\sqrt{n}}\Big)\Big)-L_D\Big(\frac{t_x}{n}\Big)\Big) = \frac{iw}{2\sqrt{n}} + \frac{1}{n} O\left( 1\right)
$$
where the $O(1)$ term is uniform in $w$. It implies that there exists a positive constant $E_{a,x}$ such that, for all $|w| \in [s_n,a\sqrt{n}]$,
\begin{equation}
\label{MAJDNXDENOM} \exp \left( n \Big(L_D\Big(\frac{1}{n}\Big(t_x+\frac{iw}{\sqrt{n}}\Big)\Big)-L_D\Big(\frac{t_x}{n}\Big)\Big)\right) \geq E_{a,x}.
\end{equation}
Then, by using \eqref{DEFCHI}, \eqref{MAJDNXNUM2} and  \eqref{MAJDNXDENOM} in \eqref{BDn}, we obtain that
\begin{eqnarray*}
|D_n(x)| &\leq & \frac{2}{t_x \sqrt{2 \pi} E_{a,x}} \int_{s_n}^{a \sqrt{n}} \frac{w}{s_n}\exp \Big( - D_{a,x} L_D^{\prime \prime}(t_x)\frac{w^2}{2} \Big) dw, \\ 
&\leq& \frac{2}{s_n t_x \sqrt{2 \pi} E_{a,x}D_{a,x} L_D^{\prime \prime}(t_x)} \exp \Big( - D_{a,x} L_D^{\prime \prime}(t_x)\frac{s_n^2}{2} \Big),
\end{eqnarray*}
which clearly means that $D_n(x)$ goes exponentially fast to zero. It remains to carry out a similar analysis for the last term $E_n(x)$. By the same lines as for inequality \eqref{BDn}, we have, up to a change of variable,
\begin{equation}
\label{MajoEn0}
 |E_n(x)| \leq 
    \frac{\chi(x,n)}{\sqrt{2 \pi n}} \int_{a \leq |w| \leq \pi n} \left \vert\exp\left(\! n \Big(L_n\Big(t_x+iw \Big) -L_n(t_x)\Big)\!\right)\right \vert dw.
\end{equation}
It follows from  \eqref{MAJDNXNUM0} together with \eqref{MAJPHID0} that
\begin{eqnarray}
\nonumber & & \hspace{-1.5cm} \frac{1}{n} \log \left \vert \exp\left(\! n \Big(L_n\Big(t_x+iw\Big) -L_n(t_x)\Big)\!\right)\!\right \vert^2\\
\label{MajoEn}  &=& \frac{1}{n} \sum_{k=1}^n \log \left(1+\frac{\exp (\frac{k t_x}{n})}{(\exp (\frac{k t_x}{n})-1)^2} 2\Big(1-\cos \Big(\frac{k w}{n}\Big)\Big)\right)\\
& & - \log \left(1+\frac{\exp (\frac{ t_x}{n})}{(\exp (\frac{ t_x}{n})-1)^2} 2\Big(1-\cos \Big(\frac{ w}{n}\Big)\Big)\right).
\nonumber
\end{eqnarray}
We already saw in Section \ref{S-SLDPDN} that for all $v \in \dR$, $2(1-\cos(v)) \leq \min(v^2,4)$. Moreover, it is well-known that for all $v\geq 0$,
$(e^v - 1)^2 \geq v^2$. Consequently, we deduce from Jensen's inequality that
\begin{eqnarray*}
& &\hspace{-1.5cm} \frac{1}{n} \sum_{k=1}^n \log \left(1+\frac{\exp (\frac{k t_x}{n})}{(\exp (\frac{k t_x}{n})-1)^2} 
2\Big(1-\cos \Big(\frac{k w}{n}\Big)\Big)\right)\\
&\leq&   \log \left(1+\frac{1}{n} \sum_{k=1}^n \frac{\exp (\frac{k t_x}{n})}{(\exp (\frac{k t_x}{n})-1)^2} 
2\Big(1-\cos \Big(\frac{k w}{n}\Big)\Big)\right), \\
&\leq& \log \left(1+\frac{e^{t_x}}{t_x^2}\frac{1}{n} \sum_{k=1}^n  \min\Big(w^2, \frac{4n^2}{k^2}\Big)\right).
\end{eqnarray*}
Therefore, we obtain from the previous inequality and a simple Maclaurin-Cauchy test that as soon as $|w|\geq a>2$,
\begin{eqnarray*}
& &\hspace{-1.5cm}\frac{1}{n} \sum_{k=1}^n \log \left(1+\frac{\exp (\frac{k t_x}{n})}{(\exp (\frac{k t_x}{n})-1)^2} 
2\Big(1-\cos \Big(\frac{k w}{n}\Big)\Big)\right)\\ &\leq& \log \left(1+\frac{e^{t_x}}{t_x^2}\int_0^1  \min\Big(w^2, \frac{4}{y^2}\Big)dy \right) \leq \log \left(1+\frac{4 e^{t_x}}{t_x^2}|w| \right).
\end{eqnarray*}
Furthermore, we have for all $v \in \dR$ such that $|v| \leq \pi$, $2(1-\cos(v)) \geq (2v/\pi)^2$. 
It implies that for all $w \in \dR$ such that
$|w| \leq \pi n$, 
\begin{eqnarray*}
& &\hspace{-1.5cm} \log \left(1+\frac{\exp (\frac{ t_x}{n})}{(\exp (\frac{ t_x}{n})-1)^2} 2\Big(1-\cos \Big(\frac{ w}{n}\Big)\Big)\right) \\
&\geq & \log \left(1+\frac{\exp (\frac{ t_x}{n})}{(\exp (\frac{ t_x}{n})-1)^2} \frac{4 w^2}{n^2\pi^2}\right)\geq  \log \left(1+\frac{4}{\pi^2t_x^2e^{t_x}} w^2\right)
\end{eqnarray*}
where the last inequality comes from the fact that for all $t>0$, $e^t-1 \leq t e^t$.
Hence, putting these two inequalities into \eqref{MajoEn} and recalling that $|w|\geq a$, we obtain that
\begin{equation}
\label{MajoEn2}
\left\vert \exp\left(\! n \Big(L_n\Big(t_x+iw\Big) -L_n(t_x)\Big)\!\right)\right\vert^2 \leq \left( \frac{1+ \frac{4 e^{t_x}}{a t_x^2}w^2}{1+\frac{4}{\pi^2t_x^2e^{t_x}} w^2} \right)^{n/2}.
\end{equation}
Hereafter, we can choose $a =2 \pi^2 e^{2t_x}>2$, leading for all $|w| \geq a$, to
\begin{equation} 
\label{FINQAX}
\left( \frac{1+ \frac{4e^{t_x}}{a t_x^2}w^2}{1+\frac{4}{\pi^2 t_x^2e^{t_x}} w^2}\right)  \leq 
\left( \frac{1+ \frac{2 }{\pi^2 t_x^2 e^{t_x}} a^2}{1+\frac{4}{\pi^2t_x^2e^{t_x}} a^2} \right).
\end{equation}
Denote by $q_{a,x}$ the upper bound in \eqref{FINQAX}. We clearly have $0<q_{a,x}<1$. Finally, we deduce from \eqref{MajoEn0} together with
\eqref{DEFCHI}, \eqref{MajoEn2} and \eqref{FINQAX} that
$$
|E_n(x)| \leq \frac{n \sqrt{2\pi n} e^{t_x}}{ t_x } \big(q_{a,x}\big)^{n/2},$$
which means that $E_n(x)$ goes exponentially fast to zero.
\demend

\vspace{-0.5cm}

\bibliographystyle{abbrv}
\bibliography{Biblio-SLDPDI}

\end{document}